\documentclass[12pt]{article}
\usepackage{amsfonts}
\usepackage{amsmath}
\usepackage{latexsym}
\usepackage{amscd}

\newtheorem{theorem}{Theorem}
\newtheorem{lemma}{Lemma}

\newtheorem{corollary}{Corollary}
\newtheorem{remark}{Remark}
\newtheorem{example}{Example}

\begin{document}
\bibliographystyle{plain}

\thispagestyle{empty}
\setcounter{page}{0}

\vspace {2cm}

{\Large G. Morvai and B. Weiss: 
 On Classifying Processes.}

\vspace {2cm}

{\Large  Appeared in:  Bernoulli  11  (2005),  no. 3, pp. 523--532.}

\vspace {2cm}

\begin{abstract}
 We prove several results concerning classifications, based on successive observations $(X_1,\dots, X_n)$ 
of an unknown stationary and ergodic 
 process, for membership in a given class of processes, such as the class of all finite order Markov chains.
\end{abstract}



\section{Introduction and Statement of  Results}

If $\cal G$ is a subclass of all stationary and ergodic binary processes then a sequence of functions 
$g_n: \{0,1\}^n\rightarrow \{YES,NO\}$ is a classification  for ${\cal G}$ in probability if 
$$
\lim_{n\to\infty} P(g_n(X_1,\dots,X_n)=YES)=1
$$
for all processes in ${\cal G}$, and 
$$
\lim_{n\to\infty} P(g_n(X_1,\dots,X_n)=NO)=1
$$
for all processes not in ${\cal G}$. 

\noindent
Similarly, 
$g_n: \{0,1\}^n\rightarrow \{YES,NO\}$ is a classification for ${\cal G}$ in a pointwise sense  if 
$$
 g_n(X_1,\dots,X_n)=YES \ \mbox{eventually almost surely}
$$
for all processes in ${\cal G}$, and 
$$
g_n(X_1,\dots,X_n)=NO \ \mbox{eventually almost surely}
$$
for all processes not in ${\cal G}$. 
Of course, if $g_n$ is a classification in a pointwise sense then it is a classification in probability 
but a classification in probability is not necessarily a classification in a pointwise sense.

For the class ${\cal M}_k$ of $k$-step mixing Markov chains of fixed order $k$, 
there is a pointwise classification of the type we have just described. 
(For mixing Markov chains see Proposition I.2.10 in Shields (1996).) 
It was carried out in detail for independent processes by Bailey (1976). 
(Actually he proved the result only for independent processes and  indicated how to 
generalize his result for 
the class of ${\cal M}_k$.) 
For the class ${\cal M}_{mix}=\bigcup_{k=0}^{\infty} {\cal M}_k$ of mixing 
Markov chains of any order, 
Bailey showed that no such classification exists. 
See Ornstein and Weiss (1990) for some further results on this kind of question. 
Our concern in this paper is with the class of finitarily Markovian processes which is defined as follows.

\noindent
Let $\{X_n\}_{n=1}^{\infty}$ be a stationary and ergodic binary time series. 
A one sided stationary time series $\{X_n\}_{n=1}^{\infty}$ can always be thought to be a two sided time series 
$\{X_n\}_{n=-\infty}^{\infty}$. For $m\le n$ let $X_m^n=(X_m,\dots,X_n)$.  

\bigskip
\noindent
{\bf Definition:} A stationary and ergodic binary time series $\{X_n\}$  is said to be finitarily Markovian if 
for  almost every $x^{-1}_{-\infty}$ there is a finite $K(x^{-1}_{-\infty})$ such that for all $i>0$ and $y^{-1}_{-i}$
if $P(X_0=1| X^{-K-1}_{-K-i}=y^{-1}_{-i},X^{-1}_{-K}=x^{-1}_{-K})>0$ then 
$$
P(X_0=1|X^{-1}_{-K}=x^{-1}_{-K})=P(X_0=1| X^{-K-1}_{-K-i}=y^{-1}_{-i},X^{-1}_{-K}=x^{-1}_{-K}). 
$$
This class includes all finite order Markov chains (mixing or not) and many other processes such as the finitarily deterministic processes of 
Kalikow, Katznelson and Weiss (1992).

\begin{example}
 First we define a 
Markov process which serves as the technical tool for  our construction.  
Let the state space $S$ be the non-negative integers. 
The transition probabilities are as follows: with probability one move from $0$ to $1$ 
and from $1$ to $2$, for all $s\ge 2$ move with equal probability $0.5$ to $0$  and $s+1$.  
This construction yields a stationary and ergodic 
Markov process $\{M_i\}$ with stationary distribution
$$
P(M_i=0)=P(M_i=1)={1\over 4}
$$
and
$$
P(M_i=j)={1\over  2^{j}} \mbox{\ \  for $j\ge 2$}.
$$
Now we define the  binary hidden Markov chain   $\{X_i\}$, which we denote as, $X_i=f(M_i)$. 
Let $f(0)=0$, $f(1)=0$, and $f(s)=1$ for all even states $s$. 
A feature of this definition of $f(\cdot)$ is that whenever  
$X_n=0,X_{n+1}=0,X_{n+2}=1$ we  know 
 that $M_n=0$ and {\it vice versa}. 
Consider the class of  processes of the above form for all possible labeling  
of the rest of the states by zero and one.
(It is easy to see  that this  class  contains Markov chains of order $\le r+1$,
 e.g. when for all $s\ge r$ $f(s)=1$ and 
processes which are not Markov of any order, e.g. when $f(2^i+1)=0$ for $i=1,2,\dots$ and 
for the rest of the yet unlabeled odd states $s$, $f(s)=1$.)  
This class is a subclass of all stationary and ergodic binary finitarily Markovian processes. 
(Clearly, the conditional probability $P(X_1=1|X^0_{-\infty})$ does not depend on values beyond the first (going backward) 
occurrence of $001$.) 
Gy\"orfi, Morvai and Yakowitz (1998) 
proved  that there is no estimator of the value
$P(X_{n+1}=1|X_1^n)$ from samples $X_1^n$ such that the error tends to zero as $n$ tends to infinity 
in the pointwise sense for this class of processes.  
\end{example}

\begin{example}
Let $\{M_n\}$ be any stationary and ergodic first order Markov chain with 
finite or countably infinite 
state space $S$. 
Let $s\in S$ be an arbitrary state with $P(M_1=s)>0$. Now let $X_n=I_{\{M_n=s\}}$. 
By Shields (1996), Chapter I.2.c.1, the binary time series $\{X_n\}$ is  stationary and ergodic. 
It is also finitarily Markovian. (Indeed, the conditional probability 
$P(X_1=1|X^0_{-\infty})$
does not depend on values  
beyond the first (going backwards) occurrence of one  in $X^0_{-\infty}$ 
which identifies the first (going backwards) occurrence of state $s$ in the 
Markov chain $\{M_n\}$. )
The resulting time series $\{X_n\}$ is not a Markov chain of any order in general. 
(Indeed, consider the Markov chain $\{M_n\}$ with state space $S=\{0,1,2\}$ and 
transition probabilities
$P(X_2=1|X_1=0)=P(X_2=2|X_1=1)=1$, $P(X_2=0|X_1=2)=P(X_2=1|X_1=2)=0.5$. This yields 
a stationary and ergodic Markov chain $\{M_n\}$, cf. Example I.2.8 in Shields (1996).  
Clearly, the resulting time series  $X_n=I_{\{M_n=0\}}$ will not be Markov of any order.
The conditional probability 
$P(X_1=0|X^0_{-\infty})$ depends on whether until the first (going backwards) 
occurrence of one 
you see
even or odd number of zeros.)    
These examples include all stationary and ergodic binary renewal processes with finite expected 
inter-arrival times, a basic class for many applications. 
(A stationary and ergodic binary renewal process is defined as 
a stationary and ergodic binary process such that the times between occurrences of 
ones  are independent and identically distributed with finite expectation, cf. 
Chapter I.2.c.1 in Shields (1996). )
\end{example}

\noindent
Our main result is that there is no classification for membership in the class 
of finitarily Markovian processes. 
As a byproduct we will also improve Bailey's result from mixing 
Markov chains to the class of Markov chains. Our results apply to both pointwise classifications and classifications in probability. 

\
\noindent
\begin{theorem} \label{thm1}
Given a sequence of functions $g_n: \{0,1\}^n\rightarrow \{YES,NO\}$ such that 
\begin{itemize}
\item for all stationary and ergodic binary Markov chains $\{X_n\}$ with arbitrary finite order
\begin{equation}
\label{assumption1}
\lim_{n\to\infty}P(g_n(X_1^n)=YES )=1
\end{equation}
\item  for all stationary and ergodic binary non finitarily Markovian processes
\begin{equation}
\label{assumption2}
\lim_{n\to\infty} P(g_n(X_1^n)=NO)=1
\end{equation}
\end{itemize}
we construct a single stationary and ergodic binary process $\{X_n\}$such that 
$$
\limsup_{n\to\infty} P(g_n(X_1^n)=YES)=1 \ \ \mbox{and } \ \ \limsup_{n\to\infty} P(g_n(X_1^n)=NO)=1.
$$
\end{theorem}

\begin{corollary}
There is no classification for the class of all stationary and ergodic binary Markov chains with arbitrary finite order, 
in a pointwise sense or in probability.
\end{corollary}

\noindent
\begin{remark}
For  motivation consider  the universal intermittent estimation problem where 
 the goal is to find stopping times $\tau_k$ 
such that one can estimate 
 $P(X_{\tau_k+1}=1|X_1^{\tau_k})$ from samples $X_1^{\tau_k}$ in the pointwise sense for all stationary and ergodic binary time series. 
Such a universal scheme was proposed in Morvai (2003). Unfortunately 
the stopping times of Morvai (2003) grow very rapidly.   
Had one  classified the 
Markov chains from non Markov chains then 
one could have improved the scheme of Morvai such that it would have remained universially 
pointwise consistent 
for all stationary and ergodic processes and particularly, 
if the process turned out to be Markov,  one could have 
estimated the conditional probability $P(X_{k+1}=1|X_1^{k})$ eventually for all $k$ that is,  
$\tau_{n+1}=\tau_n+1$ eventually. 
Indeed, if $g_n(X_1^n)$ classified the process as Markov then  one could simply use  a Markov 
order estimator  ( e.g of  Csisz\'ar and Shields (2000) ) and 
count frequencies of blocks with length equal to the order and this estimator is  
consistent in the pointwise sense 
 for Markov chains. Otherwise one could  use  the universal estimator of Morvai (2003). 
\end{remark}

\begin{corollary}
There is no classification for the class of all  stationary and ergodic binary finitarily  Markovian  processes, 
in a pointwise sense or in probability.
\end{corollary}

\noindent
\begin{remark}
Concerning the above mentioned intermittent estimation problem, 
one could  have improved 
the universal estimator of  Morvai (2003) for finitarily Markovian processes. 
Had $g_n(X_1^n)$  classified 
the process as a finitarily Markovian process  
one could use the stopping times and estimator e.g as in Morvai and Weiss (2003) 
which estimator
is not universal but it works for all finitarily Markovian processes 
and the growth of the stopping times is much more 
moderate compared to the stopping times associated with the universal estimator in Morvai (2003). 
For non finitarily Markovian processes  one could use the universal estimator of 
Morvai (2003). 
\end{remark}

\section{Proofs}

The following lemma is well known. 
\begin{lemma} \label{lemma1}
Let $\{X_n\}$ be a stationary and ergodic binary time series and  $N$  a positive integer. 
Then there is a stationary and ergodic binary Markov chain $\{Z_n\}$ of some finite  order 
$\le N$  such that  
the $N$ dimensional distributions of $\{X_n\}$ and $\{Z_n\}$ are identical. 
\end{lemma}
{\bf Proof:}
Put 
$P(Z_{N+1}=z|Z_1^N=x_1^N)=P(X_{N+1}=z|X_1^N=x_1^N)$. 
This yields a stationary and ergodic Markov chain $\{Z_n\}$ 
of some finite  order 
$\le N$ 
with 
the original marginal distribution
$P(Z_1^N=x_1^N)=P(X_1^N=x_1^N)$, that is, for $n>N$, define 
$$
P(Z_1^n=x_1^n)=P(Z_1^N=x_1^N)\prod_{i=N+1}^n P(Z_i=z_i|Z^{i-1}_{i-N}=x^{i-1}_{i-N}).
$$
Clearly $\{Z_n\}$ is a stationary Markov chain 
of some finite  order 
$\le N$ 
 since $\{X_n\}$ was stationary. 
The chain $\{Z_n\}$  can be thought of as one step Markov chain  by passing to 
$N$-tuples. The ergodicity of the $\{X_n\}$ process guarantees that this chain is 
irreducible when considered as a chain  on 
those $N$-tuples which have  positive measure under the distribution of $X_1^N$. 
The process $\{Z_n\}$ is also ergodic since stationary 
binary irreducible Markov chains of some finite order are ergodic by 
Proposition I.2.9 in Shields (1996). (Cf. also Kemeny and Snell (1960).) 
The proof of Lemma~\ref{lemma1} is complete.

\bigskip
\noindent
{\bf Definition:}
The entropy rate $H$  associated with a stationary binary time series $\{X_n\}$ is defined as 
\begin{eqnarray*}
\lefteqn{ H =
 -E \left\{ 
P(X_0=1|X^{-1}_{-\infty})\log_2 P(X_0=1|X^{-1}_{-\infty}) \right.}\\
&+&  \left. P(X_0=0|X^{-1}_{-\infty})\log_2 P(X_0=0|X^{-1}_{-\infty})
\right\}.
\end{eqnarray*}

\bigskip
\begin{lemma} \label{lemma2}
Given a stationary and ergodic binary process $\{X_n\}$, an integer $N>0$ and a real number 
$0<\delta<1$, there exists  
  a stationary and  ergodic   non finitarily Markovian process $\{Y_n\}$ such that 
 \begin{equation} \label{vardist}
\sum_{y_1^{N}\in \{0,1\}^{N} } 
|P(X_1^N=y_1^N)-
P(Y_1^N =y_1^N  )|<\delta.
 \end{equation}
\end{lemma}
{\bf Proof:}
Let $\{Z_n\}$  be a stationary and ergodic binary time series with zero entropy rate such that all finite words 
have positive probability. 
It is well known that such processes exist. For the sake of completeness we supply a proof in 
Lemma~\ref{lemmaappendix} in the Appendix.
This process is clearly not finitarily Markovian. 

\noindent
By ergodicity of the $\{X_n\}$ process, 
there  exists an $r$ and a word $w_1^r$ such that the empirical counts of all $N$ blocks from $w_1^r$ are 
$\delta/2^{N+1}$ 
close to 
the probabilities corresponding to the $\{X_n\}$ process. 

\noindent
We would like to define a process in which we alternate between  the fixed word $w_1^r$ and the 
$Z_n$'s,  $Z_1,w_1^r,Z_2,w_1^r,\dots$. If we can do this and identify uniquely the position of the 
$Z_n$'s then this process will not be finitarily Markovian. In  order to uniquely  identify  the positions 
of the $Z_n$'s we will add a synchronizing   word $u_1^m$ whose length  is very small compared to the 
length of $w_1^r$ and
which appears only where  we  place it. The fact that its length is small means that the  
finite  distributions 
will remain close to the finite distribution of the $\{X_n\}$ process. For $u_1^m$ to 
sychronize we need to 
know that when looking  across a string  like 
$Z_1,u_1^m,w_1^r,Z_2, u_1^m, w_1^r,Z_3$ the word appears only in the two locations where it is written.  

\noindent
Now choose some word $u_1^m$ with length  $m=\lceil 10\log_2 r\rceil$  such that this word $u_1^m$ 
does not appear in the word 
$w_1^r$ and 
it has no reasonable  non-trivial self overlap. More precisely, there is no non-trivial 
self overlap greater than $2/5 m$ and there is no overlap 
with $w_1^r$
 greater than $2/5m$.  The number of words with length $m$ which have  greater 
self overlap is at most $2 m 2^{3/5 m }$. 
 The number of words of length $m$  which have overlap with $w_1^r$  greater than $2/5m$ but 
not completely contained in $w_1^r$ 
 is at most $2 m 2^{3/5 m }$. The number of words with 
 length $m$ completely 
  contained in $w_1^r$ is at most $r$. Summing up the number of these possible bad  words we get 
  $$
  r+4 m 2^{3/5 m}< 2^m.
  $$  
Thus there is at least one word $u_1^m$ with the desired property. 
The word  $u_1^m$ will serve as a synchronyzing word. 

\noindent
We will define the desired $\{Y_n\}$ process in two steps. First we will define a nonstationary 
process $\{W_n\}$ as follows.
Consider $n-1=\eta(m+r+1)+\theta$, where $0\le \theta\le m+r$ and $\eta\ge 0$. 
The process $\{W_n\}$ will be obtained by inserting a fixed block $u_1^m,w_1^r$  of length 
$m+r$ between 
successive symbols of the process $\{Z_n\}$. 
Define the process $\{W_n\}$ as follows. 
Let 
\[
W_n=\left\{ \begin{array}{ll}
Z_{\eta+1} 	& \mbox{if $\theta=0$}\\
u_{\theta}	& \mbox{if $1\le \theta\le m$}\\
w_{\theta-m}	& \mbox{if $m+1\le \theta\le m+r$.}\\
\end{array}
\right. 
\]
Our assumptions on the synchronizing word  imply  that 
such  a process will not  be stationary  and to ensure stationarity  we need to randomize  over $m+r+1$.
Here is a formal description. Let $\zeta$ be distributed on $\{0,\dots, m+r\}$ uniformly. 
Let $\zeta$ be independent from $\{W_n\}$.  Define $\{Y_n\}$  as 
$Y_n=W_{n+\zeta}$. (That is, $\{Y_n\}$  is constracted from 
$\{W_n\}$ by averaging over the $m+r+1$ shifts of the $\{W_n\}$ process. )   

\noindent
The fact that $u_1^m$ was synchronyzing means that $\zeta$ is a function of the $\{Y_n\}$ process. 
Thus from $\{Y_n\}$ one recovers exactly the $\{Z_n\}$ process. 
Now $\{Y_n\}$ is a stationary and ergodic binary non  finitarily Markovian time series since $\{Z_n\}$ was such.  
To see that (\ref{vardist}) is satisfied one uses the property of $w_1^r$ and takes $r$ sufficiently large so
that the edge effects caused by $u_1^m$ are negligible.
The proof of Lemma~\ref{lemma2} is complete.

\bigskip
\noindent
{\bf Proof of Theorem \ref{thm1}:}
To construct $\{X_n\}$ we will alternately use the two lemmas to construct a sequence  of 
processes $\{Y_n^{(i)}\}$, which for odd $i$ will be a Markov chain and for even $i$ will not even be 
finitarily Markovian but the entire sequence will converge to an ergodic process $\{X_n\}$ which 
will have the required properties. Here is how this is done.  
Let $0<\epsilon_k<1$ such that $\epsilon_k\to 0$ and $0<\delta_k<1$ such that $\sum_{k=1}^{\infty} \delta_k<0.25.$
We construct our process as follows: 
Let  $\{Y_n^{(1)}\}$ 
be independent and identically distributed random variables assuming  the values $\{0,1\}$ 
with equal probabilities.  
Let $N_1>1$ be so large that 
$$
P(g_{N_1}(Y_1^{(1)},\dots,Y_{N_1}^{(1)})=YES)\ge 1-\epsilon_1
$$
and there exists  a set ${\cal U}_{N_1}\subseteq \{0,1\}^{N_1}$ such that 
$P((Y_1^{(1)},\dots,Y_{N_1}^{(1)})\in {\cal U}_{N_1})>1-\epsilon_1$ and 
$$
\max_{u_1^{N_1}, v_1^{N_1}\in {\cal U}_{N_1} }
\sum_{x\in\{0,1\} } {1\over N_1} \left|
\sum_{i=0}^{N_1-1} \left( I_{\{ u_{i+1}=x\}} -
 I_{\{v_{i+1}=x\}}\right) \right| <\epsilon_1.
$$
Assume for $k=2,\dots, i-1$ we have already defined a sequence of stationary and ergodic 
binary time series   
$\{Y_n^{(k)}\}$ and positive integers $N_k>k^2$ and sets ${\cal U}_{N_k}\subseteq\{0,1\}^{N_k}$ 
such that $P((Y_1^{(k)},\dots,Y_{N_k}^{(k)})\in {\cal U}_{N_k})>1-\epsilon_k$, 
$$
\sum_{y_1^{N_{k-1}}\in \{0,1\}^{N_{k-1}} } 
|P(Y_1^{(k-1)}=y_1,\dots,Y_{N_{k-1}}^{(k-1)}=y_{N_{k-1}} )-
P(Y_1^{(k)}=y_1,\dots,Y_{N_{k-1}}^{(k)}=y_{N_{k-1}} )|<\delta_{k-1},
$$
\begin{equation}
\label{cl1}
\max_{u_k^{N_k}, v_1^{N_k}\in {\cal U}_{N_k} }
\sum_{x_1^k\in\{0,1\}^k } {1\over N_k-k+1} \left|
\sum_{i=0}^{N_k-k} \left( I_{\{ u_{i+1}^{i+k}=x_1^k\}} -
 I_{\{v_{i+1}^{i+k}=x_1^k\}}\right) \right| <\epsilon_k,
\end{equation}
and 
\begin{itemize}
\item if $k$ is even then  $\{Y_n^{(k)}\}$ is not finitarily Markovian and 
$$
P(g_{N_k}(Y_1^{(k)},\dots,Y_{N_k}^{(k)}) = NO)\ge 1-\epsilon_k
$$
\item if $k$ is odd then  $\{Y_n^{(k)}\}$ is a Markov chain with some order and 
$$
P(g_{N_k}(Y_1^{(k)},\dots,Y_{N_k}^{(k)}) = YES)\ge 1-\epsilon_k.
$$
\end{itemize}
Now we define it for $i$. 
If $i$ is odd then 
apply Lemma \ref{lemma1} for $\{Y_n^{(i-1)}\}$ with $N_{i-1}$. Let $\{Y_n^{(i)}\}$ 
denote the resulting stationary and ergodic binary  Markov chain. 
 Now let $N_i>i^2$ be so large that 
$$ 
P(g_{N_i}(Y_1^{(i)},\dots,Y_{N_i}^{(i)}) = YES)\ge 1-\epsilon_i
$$
and there is a set ${\cal U}_{N_i}\subseteq\{0,1\}^{N_i}$ 
such that $P((Y_1^{(i)},\dots,Y_{N_i}^{(i)})\in {\cal U}_{N_i})>1-\epsilon_i$ and 
$$
\max_{u_i^{N_i}, v_1^{N_i}\in {\cal U}_{N_i} }
\sum_{x_1^i\in\{0,1\}^i } {1\over N_i-i+1} \left|
\sum_{j=0}^{N_i-i} \left( I_{\{ u_{j+1}^{j+i}=x_1^i\}} -
 I_{\{v_{j+1}^{j+i}=x_1^i\}}\right) \right| <\epsilon_i.
$$
By assumption (\ref{assumption1}) and 
the ergodicity of $\{Y^{(i)}_n\}_{n=1}^{\infty}$ there exists such an $N_i$.

\noindent
If $i$ is even then apply Lemma \ref{lemma2}  for $\{Y_n^{(i-1)}\}$ with $N_{i-1}$ and $\delta_{i-1}$. Let $\{Y_n^{(i)}\}$ 
denote the resulting non finitarily  Markovian  process. 
Now let $N_i>i^2$ be so large that 
$$ 
P(g_{N_i}(Y_1^{(i)},\dots,Y_{N_i}^{(i)}) = NO)\ge 1-\epsilon_i
$$
and there is a set ${\cal U}_{N_i}\subseteq\{0,1\}^{N_i}$ 
such that $P((Y_1^{(i)},\dots,Y_{N_i}^{(i)})\in {\cal U}_{N_i})>1-\epsilon_i$ and  
$$
\max_{u_i^{N_i}, v_1^{N_i}\in {\cal U}_{N_i} }
\sum_{x_1^i\in\{0,1\}^i } {1\over N_i-i+1} \left|
\sum_{j=0}^{N_i-i} \left( I_{\{ u_{j+1}^{j+i}=x_1^i\}} -
 I_{\{v_{j+1}^{j+i}=x_1^i\}}\right) \right| <\epsilon_i.
$$
By assumption (\ref{assumption2}) and 
the ergodicity of $\{Y^{(i)}_n\}_{n=1}^{\infty}$ there exists such an $N_i$.

\bigskip
\noindent
Now  it follows from the construction that for any $n\le N_k$ and  $k\le K$, 
$$
|P(Y_1^{(k)}=y_1,\dots, Y_{n}^{(k)}=y_n)- P(Y_1^{(K)}=y_1,\dots, Y_{n}^{(K)}=y_n)|\le  \sum_{i=k}^{\infty} \delta_i
$$
which tends to zero as $k\to \infty$. 
\bigskip
\noindent
Now define $\{X_n\}$ in the following way:  For each $n$ let 
$$
P(X_1^n=x_1^n)=\lim_{k\to\infty} P(Y_1^{(k)}=x_1,\dots, Y_{n}^{(k)}=x_n). 
$$
Clearly $\{X_n\}$ is stationary since all $\{Y^{(k)}_n\}$ were stationary.  
Since $P((X_1,\dots,X_{N_k})\in {\cal U}_{N_k})>1-\epsilon_k-\sum_{i=k}^{\infty} \delta_i$, 
$N_k>k^2$,   (\ref{cl1}) 
and Lemma~\ref{finiteblocklemma} in the Appendix, $\{X_n\}$ is also ergodic.
Now  it follows from  the construction that 
$$
|P(X_1^n=x_1^n)- P(Y_1^{(k)}=x_1,\dots, Y_{n}^{(k)}=x_n)|\le  \sum_{i=k}^{\infty} \delta_i.
$$
Thus for $k$ even, 
$$
P(g_{N_k}(X_1,\dots,X_{N_k}) = NO)\ge 1-\epsilon_k-\sum_{i=k}^{\infty} \delta_i
$$
and the right hand side tends to $1$ as $k\to\infty$. 
Similarly. when $k$ is odd, 
$$ 
P(g_{N_k}(X_1,\dots,X_{N_k}) = YES)\ge 1-\epsilon_k-\sum_{i=k}^{\infty} \delta_i
$$
and the right hand side tends to $1$ as $k\to\infty$. 
The proof of Theorem~\ref{thm1} is complete. 

{\section{Appendix}

\noindent
We present now the proofs of two fairly standard lemmas that we used before. 
\begin{lemma} \label{lemmaappendix}
There exists a stationary and ergodic time series $\{Z_n\}$ with zero entropy rate such that all finite words have 
positive probability. 
\end{lemma}
{\bf Proof:}
Let $T: [0,1]\rightarrow [0,1]$ denote the mapping $x\rightarrow x+ \alpha \ {\rm mod}\  1$ where $\alpha$ is a fixed irrational. 
Denote the Lebesgue measure on $[0,1]$ by $\mu$. For a measurable subset $A$ of $[0,1]$ let 
$\tau_A(x)=\min\{n\ge 1: T^n x \in A\}$ denote the first return time to $A$. Partition $A$ into 
$A_k=\{x\in A: \tau_A(x)=k\}$. 
Note that $T^i A_k : 0\le i<k\}$ are disjoint sets. 
We will define a particular set $A$ with the property that for all $k$ the sets $A_k$ will have positive measure.
Indeed, one can choose inductively points $\{x_n\}$ and $\delta_n>0$, $\sum_{m=n+1}^{\infty} m\delta_m<0.1 \delta_n$ 
sufficiently small so that if $I_n=[x_n-\delta_n,x_n+\delta_n]$
the $A$ defined as follows will have the required property:
$$
A=\bigcup_{n=1}^{\infty} [(I_n\bigcup T^n I_n) - [\bigcup_{m=n+1}^{\infty} \bigcup_{i=1}^{m-1} T^i I_m] ] .
$$  
It is easy to see that for all $k$, $\mu(A_k)>0$. In this case we can list all binary words with finite length, 
$\{0,1,00,01,\dots\}=\{w_1,w_2,\dots\}$, and denote by $|w_k|$ the length of $w_k$. Define a 
partition of $[0,1]$ into two sets $\{P_0,P_1\}$ by taking the $k$-th word $w_k$ in the list and assigning the first 
$|w_k|$  sets of $(T^0 A_k), (T^1 A_k), \dots, (T^{k-1} A_k)$ to $P_0$ or $P_1$  according to the symbols in $w_k$ and 
then assign to $P_0$ all remaining points in $[0,1]$. Finally define   a stationary and ergodic  binary process as 
follows: 
Choose $x$ uniformly on $[0,1]$ and set 
\[
Z_n(x)=\left\{ \begin{array}{ll}
1 	& \mbox{if $T^nx \in P_1$}\\
0	& \mbox{if $T^nx\in P_0$.}
\end{array}
\right. 
\]
It is clear that all finite words have positive probability. Furthermore it is well known that 
any process defined by an irrational rotation  
as above is stationary and ergodic and has zero entropy cf. Cornfeld {\it et al.} (1982). 
The proof of Lemma~\ref{lemmaappendix} is complete.

\bigskip
\noindent
\begin{lemma} \label{finiteblocklemma}
A binary stationary time series $\{X_n\}$ is ergodic if 
there is a sequence of positive integers $N_k>k^2$ tending to $\infty$, 
$\epsilon_k>0$ tending to zero and  a sequence of  sets ${\cal U}_{N_k}\subseteq \{0,1\}^{N_k}$
 with probability greater than $1-\epsilon_k$ such that for all $u_1^{N_k}, v_1^{N_k}\in 
{\cal U}_{N_k}$, 
\begin{equation}
\label{closeinaverage1}
\sum_{x_1^k\in\{0,1\}^k } {1\over N_k-k+1} \left|
\sum_{i=0}^{N_k-k} \left( I_{\{ u_{i+1}^{i+k}=x_1^k\}} -
 I_{\{v_{i+1}^{i+k}=x_1^k\}}\right) \right| <\epsilon_k.
\end{equation}
\end{lemma}
{\bf Proof:}
First observe that (\ref{closeinaverage1}) implies that for all $u_1^{N_k}, v_1^{N_k}\in 
{\cal U}_{N_k}$, and for all $j\le k$, 
\begin{equation}
\label{closeinaverage2}
\sum_{x_1^{j}\in\{0,1\}^{j} } 
{1\over N_k-k+1} \left|\sum_{i=0}^{N_k-k} \left( I_{\{u_{i+1}^{i+j}=x_1^{j}\}} -
 I_{\{v_{i+1}^{i+j}=x_1^{j}\}}\right) \right| <\epsilon_k.
\end{equation}
(Indeed,
\begin{eqnarray*}
{\lefteqn{
 \sum_{x_1^{j}\in\{0,1\}^{j} } 
{1\over N_k-k+1} \left|\sum_{i=0}^{N_k-k} \left( I_{\{u_{i+1}^{i+j}=x_1^{j}\}} -
 I_{\{v_{i+1}^{i+j}=x_1^{j}\}}\right) \right| }}\\
&=& 
\sum_{x_1^{j}\in\{0,1\}^{j} }
{1\over N_k-k+1} \left|\sum_{i=0}^{N_k-k} 
\sum_{x_{j+1}^{k}\in\{0,1\}^{k-j}} 
\left( I_{\{u_{i+1}^{i+k}=x_1^k\}} -
 I_{\{v_{i+1}^{i+k}=x_1^{k}\}}\right) \right|\\
&\le& 
\sum_{x_1^{j}\in\{0,1\}^{j} }
\sum_{x_{j+1}^{k}\in\{0,1\}^{k-j}} 
{1\over N_k-k+1} \left|\sum_{i=0}^{N_k-k} \left( I_{\{u_{i+1}^{i+k}=x_1^k\}} -
 I_{\{v_{i+1}^{i+k}=x_1^{k}\}}\right) \right|\\
&=&
\sum_{x_1^k\in\{0,1\}^k } {1\over N_k-k+1} \left|\sum_{i=0}^{N_k-k} 
\left( I_{\{u_{i+1}^{i+k}=x_1^k\}} -
 I_{\{v_{i+1}^{i+k}=x_1^k\}}\right) \right| 
\end{eqnarray*} 
which is, by assumption,  less than  $\epsilon_k$.)

\noindent
Now for any $M\le k$ and $u_1^{N_k},v_1^{N_k}\in {\cal U}_{N_k}$, 
\begin{eqnarray*}
\lefteqn{
\sum_{x_1^{M}\in\{0,1\}^M}  
{1\over N_k-M+1} \left|\sum_{i=0}^{N_k-M} \left( I_{\{u_{i+1}^{i+M}=x_1^{M}\}} -
 I_{\{v_{i+1}^{i+M}=x_1^M\}}\right) \right| }\\
&\le& 
\sum_{x_1^{M}\in\{0,1\}^{M} } 
{1\over N_k-k+1} \left|\sum_{i=0}^{N_k-k} \left( I_{\{u_{i+1}^{i+M}=x_1^{M}\}} -
 I_{\{v_{i+1}^{i+M}=x_1^{M}\}}\right) \right| 
{N_k-k+1\over N_k-M+1}\\
&+&{k-M\over N_k-M+1} 2^M\\
&\le& 
\epsilon_k +{k-M\over N_k-M+1} 2^M.
\end{eqnarray*}
where we used (\ref{closeinaverage2}). 
Thus for any $M\le k$ and $u_1^{N_k},v_1^{N_k}\in {\cal U}_{N_k}$, 
\begin{equation}
\label{closeinaverage3}
\sum_{x_1^{M}\in\{0,1\}^{M} } 
{1\over N_k-M+1} \left|\sum_{i=0}^{N_k-M} \left( I_{\{u_{i+1}^{i+M}=x_1^{M}\}} -
 I_{\{v_{i+1}^{i+M}=x_1^{M}\}}\right) \right| 
\le
\epsilon_k +{k-M\over N_k-M+1} 2^M.
\end{equation}

\noindent
Assume  the process $\{X_n\}$ is stationary but not ergodic.  Then for some $M$ and for some 
$a_1^M\in\{0,1\}^M$, 
$$
\lim_{n\to\infty} {1\over n} \sum_{i=0}^{n-1} I_{\{X_{i+1}^{i+M}=a_1^M\}}
$$
almost surely exists, but the limit is not a constant on any set of probability one. 
(Cf. Theorem 7.2.1 in Gray (1988).)
This means that there exist $\delta>0$ and positive integer $n_0$ such that
for all $n>n_0$ there will be sets $E_n,F_n \subseteq \{0,1\}^n$ of probability $>10 \delta$ 
such that for all $u_{1}^n\in E_n$ and $v_1^n\in F_n$,  
$$
{1\over n-M+1} \left| \sum_{i=0}^{n-M} \left(I_{\{u_{i+1}^{i+M}=a_1^M\}}-
 I_{\{v_{i+1}^{i+M}=a_1^M\}}\right) \right|>10\delta.
$$

\noindent
For $M$ and $\delta$ above choose $k$ large enough so that $M<k$,  
$\epsilon_k<0.5 \delta$, $2^M (k-M)/( N_k-M+1) <0.5 \delta$,  
and $N_k>n_0$. (Such a $k$ exists since $\epsilon_k\to 0$ and ${k\over N_k}< {1\over k}\to 0$.) 

\noindent
However this leads to a contradiction since ${\cal U}_{N_k}$ fills all but $\delta$ 
while on sets $E_{N_k}$ and 
$F_{N_k}$, which have probability at least $10 \delta$, the empirical distributions differ. 
(${\cal U}_{N_k}$ should have nonempty intersection with both $E_{N_k}$ and $F_{N_k}$ and so 
on ${\cal U}_{N_k}$ the emprical distribution should differ by $10\delta$ which contradicts 
  (\ref{closeinaverage3}) and the fact that 
$\epsilon_k+2^M (k-M)/( N_k-M+1)<\delta$. )
The proof of Lemma~\ref{finiteblocklemma} is complete. 

% -------  T E X T  E N D S  -----------------------------

\end{document}